\documentclass[preprint,10pt]{elsarticle}

\usepackage{lipsum}
\usepackage{amsmath, amssymb, amsfonts}
\usepackage{graphicx}
\usepackage{epstopdf}
\usepackage{algorithmic}
\usepackage[utf8]{inputenc}
\usepackage{float}
\usepackage{xcolor}
\usepackage{appendix}
\usepackage[ruled,vlined,linesnumbered]{algorithm2e}
\usepackage{commath}
\usepackage{mathtools}
\usepackage{verbatim}
\usepackage{url}

\journal{}

\begin{document}

\begin{frontmatter}
\title{Space-split algorithm for sensitivity analysis of discrete chaotic systems with unstable manifolds of arbitrary dimension}

\author[a1,a2]{Adam A. \'Sliwiak\corref{cor}}
\ead{asliwiak@mit.edu}
\author[a1,a2]{Qiqi Wang}
\ead{qiqi@mit.edu}
\address[a1]{Center for Computational Science and Engineering, Massachusetts Institute of Technology (MIT), 77 Massachusetts Avenue, Cambridge, MA, 02139, USA}
\address[a2]{MIT Department of Aeronautics and Astronautics}
\cortext[cor]{Corresponding author.}

\begin{abstract}
  Accurate approximations of the change of system's output and its statistics with respect to the input are highly desired in computational dynamics. Ruelle's linear response theory provides breakthrough mathematical machinery for computing the sensitivity of chaotic dynamical systems, which enables a better understanding of chaotic phenomena. In this paper, we propose an algorithm for sensitivity analysis of discrete chaos with an arbitrary number of positive Lyapunov exponents. We combine the concept of perturbation space-splitting regularizing Ruelle's original expression together with measure-based parameterization of the expanding subspace. We use these tools to rigorously derive trajectory-following recursive relations that exponentially converge, and construct a memory-efficient Monte Carlo scheme for derivatives of the output statistics. Thanks to the regularization and lack of simplifying assumptions on the behavior of the system, our method is immune to the common problems of other popular systems such as the exploding tangent solutions and unphysicality of shadowing directions. We provide a ready-to-use algorithm, analyze its complexity, and demonstrate several numerical examples of sensitivity computation of physically-inspired low-dimensional systems.  
\end{abstract}


\end{frontmatter}

\section{Introduction}\label{sec:intro}
Sensitivity analysis is a branch of computational science and engineering that investigates system's reaction to external influences. From the mathematical point of view, this response is usually measured as a derivative of some quantity of interest with respect to system's parameters. In case of chaotic systems, which describe seemingly disordered and hardly predictable phenomena, this type of analysis might be extraordinarily difficult due to the {\it butterfly effect}. The quantity of interest is usually expressed in terms a long-time average or higher-order distribution moments, rather than instantaneous values, of a certain observable $J$. In particular, the sought-after system's sensitivity with respect to a scalar parameter $s$ equals
\begin{equation}
    \label{eqn:intro-sensitivity}
    \frac{d\langle J \rangle}{ds} =: \frac{d}{ds} \left( \lim_{N\to\infty} \frac{1}{N} \sum_{k=0}^{N-1} J(x_k)\right),
\end{equation}
where $x_k$ denotes the system's state (snapshot) at time $k$. This quantity is highly desired in uncertainty quantification \cite{magri-uncertainty}, grid adaptation \cite{larsson-grid}, design optimization \cite{nielsen-optimization}, and other applications supporting advanced simulation. Practical algorithms for estimating sensitivity in the presence of chaos are essential in understanding the complex physics of turbulent flows and climate change \cite{blonigan-phdthesis, chandramoorthy-phdthesis}.  

The earliest (conventional) attempts at differentiating the statistics of an output approximate the time evolution of the solution perturbation \cite{jameson-conventional}. Indeed, these conventional methods require solving tangent/adjoint equations whose solutions represent the separation of two trajectories initiated at two different, but very close to each other, initial conditions. In the presence of chaos, however, the rate of separation is exponential and thus the tangent solutions grow in norm exponentially fast. This computational inconvenience was circumvented in ensemble methods using the concept of {\it ergodicity}. Instead of generating a long trajectory, Eyink et al. \cite{eyink-ensemble} proposed computing sensitivities of several truncated-in-time trajectories and taking the average of the partial results. While this approach does not suffer from the butterfly effect and is proven to work in different real-world chaotic systems \cite{chandramoorthy-ensemble}, large variances of the partial estimates make the ensemble methods prohibitively expensive even for medium-sized models. Yet another popular family of methods derives from the shadowing lemma \cite{pilyugin-shadowing} which, under the assumption of uniform hyperbolicity, guarantees the existence of a shadowing trajectory that lies withing a small distance to the reference solution for a long (but finite) time. The concept of shadowing was used to formulate an optimization problem to find the shadowing direction, which can be directly used to estimate the sensitivity of discrete systems \cite{wang-lssconvergence}. Least-Squares Shadowing (LSS) and its derivatives were successfully applied to various low- and large-dimensional ODE systems \cite{ni-nilss}, including a 3D turbulence model \cite{ni-jfm, blonigan-phdthesis}. In certain cases, however, shadowing solutions might be unphysical and feature dramatically different statistical behavior than the unperturbed trajectory \cite{chandramoorthy-shadowing}. Empirical evidence of the unphysicality of shadowing solutions was demonstrated in the advection-dominated regime of the Kuramoto-Sivashinsky model \cite{blonigan-ks}.


The linear response theory provides useful mathematical machinery for predicting the behavior of the system's output and its statistics in light of changing input parameters. The solution is usually represented in terms of a convolution of the linear response operator, which uses the information of the unperturbed system, and imposed forcing (input). This idea helped formulate the Fluctuation-Dissipation Theorem (FDT) \cite{kubo-fdt}, which was later utilized to construct efficient ergodic-averaging numerical schemes for sensitivities in chaotic systems \cite{gritsun-fdt, abramov-original}. Despite their computational efficiency, several FDT-based methods require specific assumptions for the statistical behavior (e.g., Gaussian ergodic measure) to reconstruct the linear response operator \cite{abramov-blended}. A breakthrough result in the field of linear response was published in \cite{ruelle-original, ruelle-corrections} by Ruelle, who derived a closed-form linear response formula for the sensitivity. The only assumption of Ruelle's theory is uniform hyperbolicity, which is quite liberal in the context of large physical systems, according to the hyperbolicity hypothesis \cite{galavotti-hypothesis, ruelle-hypothesis}. Given its universality, this expression can be translated to Monte Carlo schemes applicable to a wide class of chaotic systems. Ruelle's formula involves a sequence of integrals with respect to the Sinai-Ruelle-Bowen (SRB) measure, while the integrands include directional derivatives in the direction of system's perturbation of a certain composite function. Unfortunately, the direct evaluation of that sequence through ergodic averaging is impractical due to the exponential increase of the integrand in time (see Section 2 of \cite{sliwiak-srb} for a visualization of that problem). A natural remedy is to move the problematic differentiation operator away from the composite function through partial integration. This treatment, however, would require differentiating the SRB measure, which is absolutely continuous only on the unstable manifold. The space-split sensitivity (S3) method \cite{chandramoorthy-s3-new}, which is a novel approach of regularizing Ruelle's expression, splits the perturbation vector into two parts such that one strictly belongs to the unstable manifold allowing for partial integration, while the other one prevents the exponential increase of the remaining contribution. The authors of \cite{chandramoorthy-s3-new} showed such splitting exists and is differentiable for a general uniformly hyperbolic systems and rigorously proved the convergence of all the components of S3. In addition, various numerical examples clearly support the computational efficiency of this method applied to low-dimensional systems \cite{sliwiak-1d}.                

The main purpose of this work is to generalize the space-split algorithm, which was originally derived for systems with one-dimensional unstable manifolds. In other words, we seek a direct (non-approximate) numerical method based on the linear response theory applicable to chaotic systems with an arbitrary number of positive Lyapunov exponents (LEs). The natural appeal of S3 is that it averages recursively-generated data along a trajectory and rigorously converges to the true solution in time. The crux of the space-splitting method, therefore, are recursive and converging relations for different quantities that arose from both the perturbation vector splitting and partial integration. The SRB density gradient, which represents unstable derivative of the SRB measure, is one such a quantity. Indeed, Lebesgue integration by parts requires knowledge of the measure derivative itself \cite{sliwiak-densitygrad}. In case of one-dimensional chaos, this iterative relation directly follows from the measure preservation property involving the Frobenius-Perron operator \cite{sliwiak-1d}. If the unstable manifold is geometrically more complex, it is convenient to apply the measure-based parameterization of the unstable subspace, directly relating the directional derivative of measure with the coordinate chart, and subsequently apply the chain rule on smooth manifolds. This concept was demonstrated in the context of low-dimensional manifolds described by simple differential equations in \cite{sliwiak-densitygrad}, and later used to construct an ergodic-averaging scheme for the SRB density gradient \cite{sliwiak-srb}. The authors show that for any sufficiently smooth system, this type of manifold description facilitates the iterative computation of unstable derivatives of any differentiable quantity by solving a collection of regularized tangent equations. In this paper, we combine these two major concepts, i.e., perturbation space-splitting and measure-based parameterization, to derive a general Monte Carlo scheme for parametric derivatives of long-time averages, defined by Eq. \ref{eqn:intro-sensitivity}, in chaotic systems of arbitrary dimension.       

The introduction of this paper is followed by five sections. Section \ref{sec:ruelle} introduces Ruelle's linear response formula, highlights its main theoretical and practical aspects, and describes the concept of space-splitting (S3). Section \ref{sec:unstable} focuses on the unstable contribution of the space-splitting for discrete systems with higher-dimensional unstable manifolds. We rigorously derive general iterative formulas using the measure-based parameterization of the expanding subspace, and discuss their convergence properties. Based on our derivations/analysis, we propose the general space-split algorithm and analyze its computational complexity in Section \ref{sec:algorithm}. Subsequently, in Section \ref{sec:results}, we demonstrate convergence tests and provide several numerical examples of sensitivity computation using two low-dimensional chaotic systems. The main points of this paper are summarized in Section \ref{sec:conclusions}. 


\section{Ruelle's linear response formula: fundamental aspects and practical consequences}\label{sec:ruelle}

The purpose of this section is to review the linear response formula for discrete chaotic systems derived by Ruelle in \cite{ruelle-original,ruelle-corrections}. In addition, we introduce basic concepts from the dynamical systems theory that are critical in this work. Let us consider a parameterized discrete system,
\begin{equation}
    \label{eqn:ruelle-system}
    x_{k+1} = \varphi(x_k;s):=\varphi(x_k),
\end{equation}
governed by a diffeomorphism $\varphi:M\to M$, $M\subset\mathbb{R}^n$, $n\in \mathbb{Z}^+$, $s\in\mathbb{R}$, $k\in\mathbb{Z}$. Since we consider chaotic systems, the map $\varphi$ has at least one positive Lyapunov exponent. Chaoticity is also manifested by the exponential growth of the homogeneous tangent solutions, which represent infinitesimal perturbations of the primal solution. The rate of growth of the tangent solutions is in fact determined by LE values. We assume System \ref{eqn:ruelle-system} is {\it ergodic}, which implies the long-time average of a smooth observable $J:M\to\mathbb{R}$ can be computed in two distinct ways: 1) through averaging of the time series $\{J(x_0), J(x_1),...\}$ generated along a trajectory, or 2), as an expected value of $J$ with respect to the SRB measure $\mu$ \cite{young-srb}. Moreover, in ergodic systems, the long-time statistics do not depend on the choice of the $\mu$-typical initial condition $x_0$. The SRB measure is an invariant and physical measure that describes the statistical behavior of the system. This quantity is known to be absolutely continuous on the unstable manifold only. In the stable direction, across the expanding subspace, $\mu$ is generally singular with respect to the Lebesgue measure \cite{young-srb,sliwiak-srb}. The basic assumption of Ruelle's theory is uniform hyperbolicity. This property means that the tangent space decomposes into stable and unstable covariant subspaces everywhere on the attractor. Members of these two subspaces are solutions to the homogeneous tangent equations and their norms uniformly decay and grow, respectively, in time at an exponential rate. Recall there also exists a neutral subspace aligned with the flow in continuous-in-time systems (ODEs), which are beyond the scope of this work.            

Under the assumption of uniform hyperbolicity, Ruelle rigorously showed that 
\begin{equation}
\label{eqn:ruelle-ruelle}
    \frac{d\langle J\rangle}{ds} = \sum_{k=0}^{\infty} \int_{M}D(J\circ\varphi_k)\cdot \chi\,d\mu,
\end{equation}
where $D$ represents the differentiation (gradient) operator in phase space, while $\chi = \partial_s\varphi\circ\varphi_{-1}$ is the map perturbation vector and $\varphi_k = \varphi(\varphi_{k-1})$, $\varphi_0(x) = x$ . This result assumes the observable $J$ does not depend on the parameter $s$. If it were otherwise, the expected value of $\partial_s J$ would need to be added on the right-hand side (RHS). Although Eq. \ref{eqn:ruelle-ruelle} is rigorously true for uniformly hyperbolic diffeomorphisms, a modified version of Ruelle's expression has been empirically shown to be valid in statistical mechanics systems that violate this basic assumption \cite{lucarini-lrt}. Indeed, various high-dimensional systems describing complex physical phenomena behave as if they were uniformly hyperbolic \cite{galavotti-hypothesis,ruelle-hypothesis}. Therefore, Ruelle's expression could potentially be applied to various real-world chaotic systems encountered in science and engineering.

Although Eq. \ref{eqn:ruelle-ruelle} provides a closed-form expression for the desired quantity, its direct evaluation is impractical. One could na\"ively approximate each integral of the infinite series using the Ergodic Theorem, i.e., by taking the long-time average of the integrand along a trajectory. Note, however, that the integrand involves a derivative of the observable $J$ evaluated at $k$ time steps forward in time. Differentiating this composite function requires computing the product of the map Jacobians $D\varphi$ evaluated at $k$ consecutive time steps. Owing to the butterfly effect, this product grows exponentially with $k$ at the rate proportional to the largest LE $\lambda_1$. Assuming $J$ is sufficiently smooth, one can rigorously show that Ruelle's formula is equivalent to,
\begin{equation}
\label{eqn:ruelle-ruelle2}
    \frac{d\langle J\rangle}{ds} = \int_{M} DJ\cdot u\,d\mu,
\end{equation}
where $u_k = \partial_s x_k$, and
\begin{equation}
\label{eqn:ruelle-inhomogeneous}
    u_{k+1} = D\varphi_k\,u_k + \chi_{k+1},\;\;\;u_{0}=0.
\end{equation}
The subscript notation indicates the time index at which a given function is evaluated, e.g., $D\varphi_k = D\varphi(x_k)$ or $\chi_{k+1} = \partial_s\varphi\circ(\varphi_{-1}(x_{k+1})) = \partial_{s}\varphi(x_k)$. The computation of the sequence $\{u_0,u_1,...\}$ through the inhomogeneous tangent equation in Eq. \ref{eqn:ruelle-inhomogeneous} is ill-conditioned for the same reason. Note
\begin{equation}
\label{eqn:ruelle-note}
    u_{k} = \sum_{i=1}^{k-1} \left(\prod_{j=i}^{k-1} D\varphi_{j}\right) \chi_i,
\end{equation}
and therefore $\|u_k\|\sim \mathcal{O}(\exp(\lambda_1 k))$. The equivalence of Eq. \ref{eqn:ruelle-ruelle2}-\ref{eqn:ruelle-inhomogeneous} and Eq. \ref{eqn:ruelle-ruelle} directly follows from the chain rule applied to an ergodic system. The problem of exploding tangent solutions is absent only in the two scenarios: 1) $\lambda_1 < 0$ (non-chaotic system), or 2), the tangent solution $u$ is orthogonal to the unstable manifold. The first scenario is beyond the scope of this paper, while the second one, unfortunately, is highly unlikely.   

A natural remedy for the exploding product of Jacobians is the integration by parts applied to the original Ruelle's expression. In case of Lebesgue integrals, however, one also needs to differentiate the measure itself. Note $\chi$ does not generally belong to the unstable subspace and, therefore, the direct partial integration is prohibited. To circumvent this inconvenience, we apply the main idea of the space-split sensitivity (S3) method, proposed in \cite{chandramoorthy-s3-new}, to decompose the perturbation vector $\chi$ into two terms, $\chi = \chi_1 + \chi_2$, which we call the unstable and stable contribution, respectively. At this point, we require $\chi_1$ belongs to the $m$-dimensional unstable manifold at any typical point on the attractor. Let $q^{i}(x_k):=q^i_{k}$, $i=1,...,m$ denote an orthonormal basis of the unstable manifold at $x_k\in M$ (recall $\varphi$ has $m\geq 1$ positive LEs by assumption). Thus,  
\begin{equation}
    \label{eqn:ruelle-decomposition}
    \chi = \chi_1 + \chi_2 = \left(\sum_{i=1}^m c^i\,q^i\right) + \left(\chi - \sum_{i=1}^m c^i\,q^i\right),
\end{equation}
where $c^i$ is a scalar differentiable on unstable manifolds. This decomposition is local, i.e., both the unstable basis and scalar coefficients vary from point to point on the attractor. If we plug Eq. \ref{eqn:ruelle-decomposition} back to Ruelle's formula, we observe the first term (the one including $\chi_1$) can now be integrated by parts, regardless of the choice of $c^i$. Indeed, the unstable contribution involves $m$ directional derivatives of $J\circ\varphi^k$ in the directions indicated by the selected orthonormal basis of the unstable manifold. The second term, i.e., the stable contribution, can be rigorously re-expressed as a single Lebesgue integral, analogously to Eq. \ref{eqn:ruelle-ruelle2},
\begin{equation}
\label{eqn:ruelle-stable}
    \sum_{k=0}^{\infty} \int_{M}D(J\circ\varphi_k)\cdot \chi_2\,d\mu = \int_M DJ \cdot v \,d\mu,
\end{equation}
where $v$ is a solution to the following {\it regularized} tangent equation,
\begin{equation}
\label{eqn:ruelle-inhomogeneous-reg}
    v_{k+1} = D\varphi_k\,v_k + \chi_{k+1} - \sum_{i=1}^m c^i_{k+1}\,q^i_{k+1},\;\;\;v_{0}=0,
\end{equation}
which is derived in the same fashion as its original counterpart in Eq. \ref{eqn:ruelle-inhomogeneous}. By imposing the following set of $m$ scalar constraints,
\begin{equation}
\label{eqn:ruelle-constraint}
    v_{k+1}\cdot q_{k+1}^{i} = 0,\;i=1,...,m,
\end{equation}
we can define the scalar coefficients $c^i$ that enforce the tangent solution $v$ to be orthogonal to the unstable manifold everywhere along a trajectory. This guarantees the norm of $v$ does not increase exponentially in time. Therefore, by combining Eq. \ref{eqn:ruelle-inhomogeneous-reg}-\ref{eqn:ruelle-constraint}, we obtain a linear system with $n+m$ equations and the same number of unknowns ($n$ components of $v$ and $m$ scalars $c^i$). Using the fact $q^i\cdot q^j = 1$ if $i=j$ and $q^i\cdot q^j = 0$ otherwise, the orthogonality constraint can be enforced by setting
\begin{equation}
\label{eqn:ruelle-c}
    c_{k+1}^i = q^i_{k+1}\cdot\left(D\varphi_k\,v_k + \chi_{k+1}\right),\;i=1,...,m.
\end{equation}

Due to the ergodicity and uniform hyperbolicity of $\varphi$, the RHS of Eq. \ref{eqn:ruelle-stable} can be approximated as a finite-time average of $DJ\cdot v$. Moreover, assuming the integrand is H\"older continuous, this approximation rigorously converges to the true solution as $\mathcal{O}(1/\sqrt{N})$, i.e.,
\begin{equation}
\label{eqn:ruelle-stable-convergence}
    \left|\int_M DJ \cdot v \,d\mu - \frac{1}{N}\sum_{k=0}^{N-1}DJ_k\cdot v_k\right|\lessapprox \frac{C}{\sqrt{N}}.
\end{equation}
The reader is referred to \cite{chandramoorthy-s3-new} (Proposition 8.1) for a proof of the preceding statement. We use the symbol ``$\lessapprox$" to acknowledge this estimate is approximate. Indeed, the rigorous result also involves a factor proportional to $\log\log N$ in the numerator, which is approximately constant if $N$ is large. 

Therefore, the computation of the stable contribution requires solving a constrained tangent equation (Eq. \ref{eqn:ruelle-inhomogeneous-reg}-\ref{eqn:ruelle-constraint}) and taking the time average of $DJ\cdot v$, which we compute as we travel along a trajectory initiated at a $\mu$-typical point. Note the dimension of the unstable manifold has little impact on the complexity of the algorithm and thus on its total computational cost. To enforce all the orthogonality constrains, one naturally needs to pre-compute the basis of the unstable manifold everywhere along the trajectory. We postpone the discussion on numerical procedures for approximating the unstable basis vectors until Section \ref{sec:unstable}, as they are essential in the unstable contribution.  

In conclusion, the perturbation vector space-splitting described by Eq. \ref{eqn:ruelle-decomposition} is absolutely critical is regularizing Ruelle's linear response expression for a chaotic system's sensitivities. Indeed, it allows to partially integrate the term involving the unstable component of $\chi$, and apply the conventional tangent equation approach together with recursive orthogonalization to approximate the remaining (stable) term. The algorithm for approximating the latter is largely agnostic to the dimensionality of the unstable manifold. This not the case for the unstable contribution, which is the main focus of the following section.

\section{Computation of the unstable contribution}\label{sec:unstable}

To approximate the desired sensitivity, $d\langle J\rangle/ds$, one needs to sum up two contributions, the unstable and stable terms, as defined in the previous section (Eq. \ref{eqn:ruelle-decomposition}). Since the former contains a component of the perturbation vector that is a member of the unstable manifold, we are allowed to apply partial integration to move the differentiation operator away from the troublesome composite function. However, we can do so only on unstable manifolds, because the SRB measure $\mu$ is generally non-smooth on $M$. Therefore, an extra step involving measure disintegration is required before applying the integration by parts. Let us consider a measurable partition $U$ of $M$ aligned with the geometry of unstable manifolds. Thus, for any Borel subset $B\subset M$,
\begin{equation}
    \label{eqn:unstable-disintegration}
    \mu(B) = \int_{M/U} \tilde{\mu}_x(B\cap U_x)\,d\hat{\mu}(x),
\end{equation}
where $\hat{\mu}$ denotes the quotient measure defined by the partition $U$, while $\tilde{\mu}_x$ represents the SRB measure conditioned on the unstable manifold $U_x$ that contains $x\in M$. Intuitively, Eq. \ref{eqn:unstable-disintegration} means that the measure of $B$ can be computed by summing conditional measures of local intersections weighted by the likelihood of each partition member. Using Eq. \ref{eqn:unstable-disintegration}, the unstable contribution can be expressed as follows,
\begin{equation}
    \label{eqn:unstable-unstable}
    \sum_{k=0}^{\infty} \int_{M}D(J\circ\varphi_k)\cdot \chi_1\,d\mu = \sum_{k=0}^{\infty} \sum_{i=1}^m \int_M c^i\,\partial_{q^i}(J\circ\varphi_k)\,d\mu,
\end{equation}
where $\partial_{q^i}(\cdot):=D(\cdot)\cdot q^i$ is a short-hand notation for the directional derivative in the direction of the $i$-th basis vector. Within each Lebesgue integral of the above double sum, we apply measure disintegration and then integrate by parts on $U_x$,
\begin{equation}
    \label{eqn:unstable-partial-integration}
    \begin{split}
    &\int_M c^i\,\partial_{q^i}(J\circ\varphi_k)\,d\mu =  \\& \int_{M/U}\int_{U_x} c^i\,\partial_{q^i}(J\circ\varphi_k)\, d\tilde{\mu}_x\,d\hat{\mu}(x)= 
    \int_{M/U}\int_{U_x} c^i\,\partial_{q^i}(J\circ\varphi_k)\,\tilde{\rho}_x\,d\tilde{\omega}_x\,d\hat{\mu}(x)= \\&
    -\int_{M/U}\int_{U_x} J\circ\varphi_k\,\left(\partial_{q^i} c^i + c^i\frac{\partial_{q^i}\tilde{\rho}_x}{\tilde{\rho}_x}\right)\,d\tilde{\mu}_x\,d\hat{\mu}(x) + \mathrm{B.T.},
    \end{split}
\end{equation}
where $\tilde{\rho}_x$ and $\omega_x$ respectively represent the density of the conditional measure and the natural volume form, both defined on $U_x$. The second term on the RHS of Eq. \ref{eqn:unstable-partial-integration} represents the boundary term, denoted by $\mathrm{B.T.}$, which can be expressed as the divergence of a smooth field on unstable manifolds. In uniformly hyperbolic systems, this term rigorously vanishes according to Theorem 3.1(b) of \cite{ruelle-original}. The reader is also referred to \cite{sliwiak-srb} for a more intuitive explanation of this counter-intuitive cancellation. We eventually obtain a new integral that involves two quantities, $b^{(i,i)}$ and $g^i$, defined as follows,
\begin{equation}
    \label{eqn:unstable-new-quants}
    b^{(i,j)} := \partial_{q^j} c^i,\;g^i := \frac{\partial_{q^i}\tilde{\rho}_x^s}{\tilde{\rho}_x^s} = \partial_{q^i}\log \tilde{\rho}_x^s.
\end{equation}
The computation of these two quantities is the actual price for the regularization of the original Lebesgue integrals. The latter is known in the literature as the SRB density gradient \cite{sliwiak-1d, sliwiak-differentiability, chandramoorthy-s3-new}. It reflects a relative measure change along an unstable manifold and its value is thus independent from its corresponding quotient measure. An efficient trajectory-driven algorithm for the computation of $g$ has been proposed the authors in \cite{sliwiak-srb}. We intend to utilize the measure-based parameterization of unstable manifolds proposed in that work to derive an efficient algorithm for computing the unstable contribution. In other words, we apply the machinery of iterative step-by-step computation of directional derivatives appearing in the RHS of Eq. \ref{eqn:unstable-partial-integration}. For this purpose, let us consider a family of smooth charts $x_{k}(\xi):[0,1]^m\to U_k$, where $U_k$ is a partition member that is crossed by the trajectory at time step $k$. In particular, $U_k$ is an $m$-dimensional unstable manifold such that $x_k\in U_k$. The measure-based parameterization is defined such that the SRB measure of any Borel subset $V\in [0,1]^m$ satisfying $x_{k}(V)=B_k\subset U_k$ is related to the corresponding SRB density through
\begin{equation}
    \label{eqn:unstable-parameterization}
    \tilde{\mu}_k(V) = \int_{B_k} \tilde{\rho}_k\;d\omega_k.
\end{equation}
Notice we replaced $x$ with $k$ in the subscript of the local/conditional quantities. The integer $k$ indicates the time step and thus, for example, $\tilde{\rho}_k$ denotes the conditional SRB density defined on an unstable manifold containing $x_k$. Indeed, for a given trajectory (initial condition), $k$ uniquely determines a point on the attractor. The major benefit of this type of description is a straightforward relation between the parametric gradient of $x_k(\xi)$, denoted by $\nabla_{\xi} x_k$, and the conditional SRB density $\tilde{\rho}_k$. In particular,
\begin{equation}
    \label{eqn:unstable-measure-convervation}
    \tilde{\rho}_k\,|\det R(x_k)| = 1
\end{equation}
for any $\xi\in [0,1]^m$, where $R$ is an $m\times m$ invertible matrix obtained through the QR factorization (orthonormalization) of the chart gradient,
\begin{equation}
    \label{eqn:unstable-qr}
\nabla_{\xi} x_k = Q(x_k)\,R(x_k),
\end{equation}
where $Q^T\,Q=I$ and $R$ is an upper-triangular matrix containing projections of the columns of $\nabla_{\xi} x_k$ onto its orthonormal basis stored in the $Q$ matrix. Eq. \ref{eqn:unstable-measure-convervation} is a general representation of the measure conservation of a nonlinear transformation from a uniform (constant) to non-uniform distribution. Taking the directional derivative of Eq. \ref{eqn:unstable-measure-convervation} in the direction of the $i$-th basis vector $q^i_k$, which is stored in the $i$-th column of $Q_k$, one can derive a closed-form expression for $g^i_k$ \cite{sliwiak-densitygrad},
\begin{equation}
    \label{eqn:unstable-g-full}
    g^i_k = -\frac{\mathrm{tr}\left(Q^T(x_k) \;\partial_{\xi^{(i)}}\nabla_{\xi}x_k\;R^{-1}(x_k)\right)}{\|\partial_{\xi^{(i)}}x_k\|},
\end{equation}
which is valid for all $\xi\in [0,1]^m$. The operator $\partial_{\xi^{(i)}}$ indicates differentiation with respect to the indicated paramater, and thus $\partial_{\xi^{(i)}}\nabla_{\xi}x_k$ is a matrix, whose columns contain second-order parametric derivatives. To reference specific components of an array, we introduce the round-bracket notation in the superscript. For example, the $i$-th component of some vector $v$ will be denoted as $v^{(i)}$, while $A^{(ij)}$ represents the entry from the $i$-th row and $j$-th column of a matrix $A$. We shall occasionally use the colon notation to reference all components of an array, for example, $Q^{(:j)} = q^j$. 

The crux of the density gradient computation, as explained in \cite{sliwiak-densitygrad}, relies on a recursive computation of the first- and second-order derivatives of the chart. These recursive formulas are simply derived by taking derivatives of the original system (Eq. \ref{eqn:ruelle-system}) and applying the chain rule. Note, however, that the na\"ive computation of the chart gradient is ill-conditioned, because $\nabla_{\xi_0}x_k = (\prod_{i=0}^{k-1} D\varphi_i)\nabla_{\xi_0}x_0$ grows in norm exponentially fast as discussed in Section \ref{sec:ruelle}, where $\xi_0$ represents the original (i.e., chosen at $k=0$) parametric coordinate system. In their recent work \cite{sliwiak-srb}, the authors proposed a step-by-step orthonormalization of the chart gradient through a recursive update of the coordinate system using the following linear transformation,
\begin{equation}
    \label{eqn:unstable-coordinate-trans}
    \xi_{k+1} = R(x_{k+1}(0))\,\xi_k.
\end{equation}
Applying this coordinate change in a step-by-step manner, we ensure the parametric gradient computed with respect to the new coordinates is orthogonal at the origin ($\xi = 0$). In practice, this requires performing the QR factorization every time step, where $Q$ contains the orthogonal basis, while $R$ is used to transform coordinates. A useful property of the formula for $g$ (Eq. \ref{eqn:unstable-g-full}) is its immunity to any linear coordinate transformation \cite{sliwiak-densitygrad}. In other words, Eq. \ref{eqn:unstable-g-full} is still valid in the locally orthogonalized system, which means that this formula can be dramatically simplified to 
\begin{equation}
    \label{eqn:unstable-g-simple}
    g^i_k = - q_k^j \cdot \partial_{\xi_k^{(i)}}\partial_{\xi_k^{(j)}}x_k:= - q_k^j \cdot a_k^{i,j}
\end{equation}
only at $\xi = 0$ \cite{sliwiak-srb}, where the repeated indices imply summation (per Einstein's convention), while $a$ satisfies the following recursion 
\begin{equation}
    \label{eqn:unstable-a}
    a^{i,j}_{k+1} = \left(D^2\varphi_k(q_k^p, q_k^q) + D\varphi_k\,a_k^{p,q} \right)(R_{k+1}^{-1})^{(pi)}\,(R_{k+1}^{-1})^{(qj)}.
\end{equation}
The product $D^2\varphi(a,b)$ represents the contraction of the Hessian of $\varphi$ against two vectors, $a$ and $b$. This operation outputs a vector whose $i$-th component equals $(D^2\varphi(a,b))^{(i)} = \partial_{x^{(p)}}\partial_{x^{(q)}}\varphi^{(i)}\,a^{(p)}\,b^{(q)}$. Note the choice of $\xi = 0$ does not restrict our algorithm to a certain trajectory; one can freely stretch/shrink the feasible space of $\xi$ such that the preimage of the initial state $x_0$ is $\xi_0 = 0$.

We shall now analyze the convergence of the recursive algorithm for the SRB density gradient. Based on the above description, the iterative computation of the basis matrix $Q$ involves two steps, i.e.,  left-multiplying $Q$ by the Jacobian matrix followed by QR factorization of the obtained matrix product. This implies that the basis matrix at the $k$-th time step equals
\begin{equation}
    \label{eqn:unstable-q-conv}
    Q_{k} = D\varphi_{k-1}\,...\,D\varphi_{0}\,Q_0\,R_{1}^{-1}\,...\,R_{k}^{-1}.
\end{equation}

If one replaces $Q_0$ with any arbitrary matrix that is bounded in norm, then the process described by Eq. \ref{eqn:unstable-q-conv} is guaranteed to converge at an exponential rate if $\varphi$ is a uniformly hyperbolic diffeomorphism \cite{kuptsov-lyapunov}. It means that in ideally chaotic systems one can generate the basis vectors of unstable manifolds (a.k.a. the backward Lyapunov vectors) by running a trajectory-driven iteration described above. In such systems, the product of the inverses of subsequent $R$ matrices decays in norm at an exponential rate (or faster) as $k$ increases \cite{ershov-lyapunov}, i.e.,
\begin{equation}
    \label{eqn:unstable-rconv}
    \|R_1^{-1}\,R_2^{-1}...R_k^{-1}\|\leq \exp(-ck),\;c>0.
\end{equation}
These matrices counterbalance the exploding product of Jacobians along typical trajectories. 

Given this remarkable behavior, we conclude that the iterative process for $a$ must also converge. To see that, let us consider a difference between two approximations of $a$ along a single trajectory assuming the basis vectors are the same in both the iterations, labelled as $1$ and $2$,
\begin{equation}
    \label{eqn:unstable-a-diff}
    a_{k+1,1}^{(i,j)} - a_{k+1,2}^{(i,j)} := \delta a_{k+1}^{(i,j)}= D\varphi_k\,\delta a_{k}^{p,q}\,(R_{k+1}^{-1})^{(pi)}\,(R_{k+1}^{-1})^{(qj)}.
\end{equation}
We observe Eq. \ref{eqn:unstable-a-diff} describes the evolution of the differences of acceleration vectors along a trajectory. This equation implies that the differences are recursively left-multiplied by the map Jacobian combined with a double contraction against the $R$ matrix. Note the RHS of Eq. \ref{eqn:unstable-a-diff} can viewed as a two-step algebraic process. In the first step, one computes $n$ matrix products $R_{k+1}^T\,\delta A_k^{i}\,R_{k+1}$, where $(\delta A^{s}_k)^{(pq)}$ is an $m\times m$ matrix that contains the $s$-th components of $\delta a^{p,q}_k$, $p,q=1,...,m$. Subsequently, all new $m^2$ $n$-dimensional vectors are left-multiplied by the same Jacobian matrix. Note these two algebraic operations are commutative, which means that we are allowed to take the initial (bounded) differences $\delta a^{p,q}_0$, $p,q=1,...,m$, left-multiply them by a product of $k$ Jacobians and then recursively compute the double contractions against $k$ inverses of $R$. Note also that the double contraction operation can be split into two single ones, which further implies we could, for example, recursively left-multiply the difference vector by a Jacobian with a single contraction, which is equivalent to replacing $(R_{k+1}^{-1})^{pi}$ with a Kronecker delta $\delta^{pi}$. To obtain the true solution at time step $k+1$, the obtained vectors would need to be recursively contracted against $k+1$ inverses of $R$ once more. The purpose of this discussion is to argue that if we replace one contraction with an identity operation in Eq. \ref{eqn:unstable-a-diff}, we effectively obtain a recursion equivalent to the one in Eq. \ref{eqn:unstable-q-conv}, which produces vectors with norms of the order $\mathcal{O}(1)$. The second contraction appearing in the original version of Eq. \ref{eqn:unstable-q-conv} means that these vectors are left-multiplied by the product $R_{k+1}^{-1}...R_{1}^{-1}$ whose induced norm uniformly approaches 0 at an exponential rate. Therefore, if the iteration defined by Eq. \ref{eqn:unstable-q-conv} exponentially converges to the true solution regardless of the choice of $Q_0$ (which is the case in uniformly hyperbolic systems), the recursion for $a$ (Eq. \ref{eqn:unstable-a}) also exponentially converges to its true value. The remarkable implication is that the iterative algorithm for the SRB density gradient $g$ does not depend on the initial guess and its true value can be obtained after a moderately small number of iterations.

The remaining part of this section focuses on recursive computation of $b$, which is the final term required to evaluate the RHS of the regularized unstable contribution in Eq. \ref{eqn:unstable-partial-integration}. Recall $b$ equals a parametric derivative of the scalars appearing in the constrained tangent equation. These scalars are directly computed using Eq. \ref{eqn:ruelle-c}. Recall also we describe the unstable manifold using a smooth chart $x_k(\xi_k)$ with a linearly rescaled coordinated system ensuring the orthogonality of its gradient at $\xi_k = 0$ as introduced above. Thus, by differenting Eq. \ref{eqn:ruelle-c}, one can obtain an explicit formula for $b$,   
\begin{equation}
    \label{eqn:unstable-b}
    \begin{split}
    &b^{i,j}_{k+1} = \partial_{q_{k+1}^j}c^i_{k+1} \stackrel{\|\partial_{\xi_{k+1}^j}x_{k+1}(0)\|=1}{=\joinrel=\joinrel=\joinrel=\joinrel=\joinrel=} \partial_{\xi_{k+1}^{(j)}}c^i_{k+1} =  \\ &
    \partial_{\xi_{k+1}^{(j)}}q_{k+1}^i\cdot f_{k} +
    q_{k+1}^i\cdot \partial_{\xi_{k+1}^{(j)}}f_{k} :=
    p^{i,j}_{k+1}\cdot f_{k} + q_{k+1}^i\cdot \partial_{\xi_{k+1}^{(j)}}f_{k},
    \end{split}
\end{equation}
where $f_k:= D\varphi_k\,v_k + \chi_{k+1} = D\varphi_k\,v_k + \partial_s\varphi_k$. While the recipe for $f$ and $q$ has already been discussed, we still require two more quantities, $p$ and parametric derivative of $f$, in order to complete the algorithm. We first focus on $p$, which equals the directional derivative of a backward Lyapunov vector at the origin of the updated coordinate system. In general, however, $p^{i,j}$ does not equal $a^{i,j}$. The latter is defined as the second parametric derivative of the chart evaluated at the origin. Our new quantity $p$, on the other hand, is defined as the parametric derivative of $q$ also evaluated at the origin. Therefore, to relate these two quantities, one also needs to differentiate the rescaling factor represented by the $R$ matrix. This relationship can be found by differentiating Eq. \ref{eqn:unstable-qr} with respect to the $i$-th chart coordinate,  
\begin{equation}
    \label{eqn:unstable-qr-diff1}
    \partial_{\xi_{k+1}^{(i)}}(\nabla_{\xi_{k+1}}x_{k+1}) = (\partial_{\xi_{k+1}^{(i)}}Q_{k+1})\,R_{k+1} + Q_{k+1}\,(\partial_{\xi_{k+1}^{(i)}}R_{k+1}),
\end{equation}
which implies that
\begin{equation}
    \label{eqn:unstable-qr-diff2}
    (\partial_{\xi_{k+1}^{(i)}}R_{k+1})\,R_{k+1}^{-1} = Q_{k+1}^T\,\partial_{\xi_{k+1}^{(i)}}(\nabla_{\xi_{k+1}}x_{k+1})\,R_{k+1}^{-1} - Q_{k+1}^T\,(\partial_{\xi_{k+1}^{(i)}}Q_{k+1}) ,
\end{equation}
for any $\xi_{k+1}\in[0,1]^m$. At the origin of the orthonormalized coordinate system, however, the $R$ matrix equals the identity by construction and thus
\begin{equation}
    \label{eqn:unstable-R-deriv}
    \partial_{\xi_{k+1}^{(i)}}R_{k+1} = Q_{k+1}^T\,A_{k+1}^i - Q_{k+1}^T\,P_{k+1}^i,
\end{equation}
where $A_{k+1}^i$ and $P_{k+1}^i$ respectively contain second parametric derivatives of the chart and first parametric derivatives of the basis vector, both evaluated at the origin. Note that $(A_{k+1}^i)^{(:j)} := a^{j,i}_{k+1} = a^{i,j}_{k+1}$ assuming $x(\xi)$ is sufficiently smooth and, analogously, $(P_{k+1}^i)^{(:j)} := p^{j,i}_{k+1}$. Although Eq. \ref{eqn:unstable-R-deriv} provides an explicit relation between $A$ and $P$, we still need more information to compute the latter as the parametric derivative of $R$ is unknown. The missing puzzle piece is hidden is the structure of the matrices appearing in Eq. \ref{eqn:unstable-R-deriv}. Indeed, the LHS of that equation is always upper-triangular by construction, while the second term on the RHS must be skew-symmetric (differentiate $Q_{k+1}^T\,Q_{k+1}=I$ to see it). Therefore, we infer that
\begin{equation}
    \label{eqn:unstable-R-deriv-comp}
    (\partial_{\xi_{k+1}^l}R_{k+1})^{(ij)} = \begin{cases} q_{k+1}^i\,\cdot a_{k+1}^{j,l} & \text{if } i = j, \\ q_{k+1}^i\cdot a_{k+1}^{j,l} + q_{k+1}^{j}\cdot a_{k+1}^{i,l} & \text{if } i < j \\
    0 & \text{otherwise}.
    \end{cases}
\end{equation}

We now combine Eq. \ref{eqn:unstable-R-deriv-comp} and Eq. \ref{eqn:unstable-qr-diff1} to infer an explicit expression for derivatives of backward Lypaunov vectors at the origin,
\begin{equation}
    \label{eqn:unstable-p}
    p_{k+1}^{i,j} = a_{k+1}^{i,j} - q_{k+1}^{l}(\partial_{\xi_{k+1}^j}R_{k+1})^{(li)}.
\end{equation}
Note the computation of $p$ requires only the knowledge of $a$ and $Q$, both of which are integral components of the algorithm for the SRB density gradient $g$ \cite{sliwiak-srb}. Therefore, the procedure for $g$ extended by the two above equations, Eq. \ref{eqn:unstable-R-deriv-comp} and Eq. \ref{eqn:unstable-p}, enables recursive computation of $p$ along a typical trajectory. If the procedure for $g$ converges exponentially fast, as argued above, the same is true of its extended version. Eq. \ref{eqn:unstable-R-deriv-comp} clearly indicates that, in general, $p^{i,j} \neq p^{j,i}$ if $i \neq j$. The lack of symmetry requires us to compute all $m^2$ different $p$ vectors to advance the full algorithm in time, which will be evident at the end of this section.

The final task in the derivation of the full algorithm is to apply the chain rule in the second term of the RHS of Eq. \ref{eqn:unstable-b}. Notice that at $\xi_{k+1} = \xi_k = 0$, one can directly change variables of the differentiation because $R_{k+1}^{-1} = \partial \xi_{k}/ \partial \xi_{k+1}$, which implies that
\begin{equation}
    \label{eqn:unstable-tranformation}
    \nabla_{\xi_{k+1}}f_{k} = \nabla_{\xi_{k}}f_{k}\,R_{k+1}^{-1},
\end{equation}
where the $i$-th column of $\nabla_{\xi_{k}}f_{k}$ can be expanded as follows,
\begin{equation}
    \label{eqn:unstable-f-chainrule}
    \partial_{\xi_k^{(i)}}\,f_k = D^2\varphi_k(v_k,q_k^i) + D\varphi_k\,w_k^i + D\partial_s\varphi_k\,q_k^i.
\end{equation}
The matrix $D\partial_s\varphi_k$ represents the Jacobian of the map differentiated with respect to the scalar $s$ and evaluated at time $k$. The new quantity, $w_k^i$, is defined as $w_k^i:=\partial_{\xi_k^{(i)}}\,v_k$ and is recursively computed in the following way,
\begin{equation}
    \label{eqn:unstable-w}
    w_{k+1}^i = (\nabla_{\xi_{k+1}}f_k)^{(:i)} - b_{k+1}^{l,i}\,q_{k+1}^l + c_{k+1}^l\,p_{k+1}^{l,i}.
\end{equation}
This formula is obtained through parametric differentiation of Eq. \ref{eqn:ruelle-inhomogeneous-reg}. We now observe the entire set of $m^2$ scalars $b^{i,j}$ and $m^2$ vectors $p^{i,j}$ are necessary in order to advance the iteration for $w^i$. While $b$ appears in the recursion for $w$ and vice versa, there is no need to construct large linear systems to find both the quantities. Indeed, Eq. \ref{eqn:unstable-b}, \ref{eqn:unstable-tranformation}, \ref{eqn:unstable-f-chainrule} indicate that in order to find $b$ at time $k+1$, we need all vectors $w$ at the previous time $k$. Therefore, in our algorithm, we can sequentially compute all vectors/scalars in the following order: $a$, $p$, $b$ and $w$, at every point along a trajectory. 

We already discussed the convergence of the iterations for $Q$, $a$, and $p$. Our final task is the convergence analysis of the recursion for $w$ (Eq. \ref{eqn:unstable-w}). Let $\delta W_k:= W_{k,1}-W_{k,2}$ be the difference of two matrices containing all vectors $w$ in their columns such that $(W_{k,1})^{(:i)}:=w_{k,1}^i$ and $(W_{k,2})^{(:i)}:=w_{k,2}^i$, while the labels 1 and 2 represent two different (randomly chosen) initial conditions for the recursion of $w$. Therefore, $\|\delta W_0\|\neq 0$ in general. Using this notation and combining Eq. \ref{eqn:unstable-w} and Eq. \ref{eqn:unstable-b}, we derive the following iteration for the difference matrix,
\begin{equation}
    \label{eqn:unstable-w-conv}
    \delta W_{k+1} = \left(I - Q_{k+1}\,Q_{k+1}^T\right)D\varphi_k\,\delta W_k\,R_{k+1}^{-1}.
\end{equation}
Note the difference matrix is left-multiplied by another matrix that is orthogonal to the unstable manifold, because $Q_{k+1}^T\left(I - Q_{k+1}\,Q_{k+1}^T\right)(\cdot) = 0$. Therefore, the recursive application of the left-hand side operator $\left(I - Q_{k+1}\,Q_{k+1}^T\right)D\varphi_k$ does not lead to the exponential growth of the resulting product in time. In fact, this product alone approaches 0 in norm exponentially fast. To see it, let $v = C_u v + C_s v$ be a generic bounded-in-norm vector in the tangent space, while $C_u v$ and $C_s v$ are its components belonging to the unstable and stable manifolds. Thus, $(I-QQ^T)C_u v = 0$ and $(I-QQ^T)C_s v = C_s v$. In addition, uniform hyperbolicity guarantees that the product $(\prod_{k=0}^N D\varphi_k) C_s v$ strictly belongs to the stable subspace (covariance property), while its norm is upperbounded by $C\lambda^N\|v\|$ with $C>0$ and $\lambda\in(0,1)$ (uniform decay property).
Note also that the initial difference $\delta W_0$ is bombarded by the product of the inverses of $R$ matrices, which also decays exponentially in norm with $k$, per our discussion above. The ultimate conclusion of this analysis is that all the recursions derived in this section do not depend on initial conditions and their respective solutions converge to their true values exponentially fast.  

Having the collection of converging iterative expressions for different quantities arising in the regularized version of the unstable contribution, the final step is to take the time average of the series generated along a typical trajectory. In particular, assuming the system is ergodic and combining Eq. \ref{eqn:unstable-unstable}-\ref{eqn:unstable-partial-integration}, we approximate the unstable contribution through the following triple sum,
\begin{equation}
    \label{eqn:unstable-approximation}
    \sum_{t=0}^{\infty} \int_{M}D(J\circ\varphi^s_k)\cdot \chi_1\,d\mu \approx \sum_{t=0}^K\sum_{k=0}^N\sum_{i=1}^m J_{k+t}\,\left(b_k^{(i,i)}+c_k^i\,g_k^i\right),
\end{equation}
where $T,K$ are some sufficiently large positive integers. Assuming all the quantities appearing in the above integral are H\"older continuous, the law of iterated logarithm applies and the truncated series approximating ergodic averages converge as $\mathcal{O}(\log\log N/\sqrt{N})$. Moreover, under the same assumption, the authors of \cite{chandramoorthy-s3-new} rigorously prove that a truncated series of ergodic averages in the form of the RHS of Eq. \ref{eqn:unstable-approximation} converges to the true solution as $N\to\infty$ followed by $K\to\infty$ (the double limit must be in that order). In the same work, the authors estimate the upperbound of the truncation error in terms of $N$ and $K$,
\begin{equation}
    \label{eqn:unstable-approximation-error}
    \begin{split}
    &\left|\sum_{t=0}^{\infty}\int_{M}D(J\circ\varphi^s_t)\cdot \chi_1\,d\mu-\sum_{t=0}^K\sum_{k=0}^N\sum_{i=1}^m J_{k+t}\,\left(b_k^{(i,i)}+c_k^i\,g_k^i\right)\right|\\
    &\lessapprox C_1\frac{K}{\sqrt{N}}+C_2\exp(-C_3\,K)
    \end{split}
\end{equation}
for some positive real constants $C_1$, $C_2$ and $C_3$. It implies that, for a fixed value of K (i.e., number of terms in the truncated series), our recursive method based on the perturbation vector splitting approximately behaves as a typical Monte Carlo algorithm. The bias associated with the truncation of the infinite series decays exponentially with $K$ only if $N\to\infty$. The summary of the entire algorithm, analysis of its computational complexity, and demonstration of numerical examples are presented in the following two sections.


\section{Space-split algorithm for chaotic maps}\label{sec:algorithm}

We synthesize all derivations and analysis presented in Section \ref{sec:ruelle} and \ref{sec:unstable}, and construct an algorithm for sensitivity computation of chaotic dynamical systems with an arbitrary number of $n$ degrees of freedom and positive Lyapunov exponents $m$. Algorithm \ref{alg:alg1} is a summary of the space-split procedure in the form of a pseudocode. 
\begin{algorithm}\label{alg:alg1}
\SetAlgoLined
\SetKwInOut{Input}{Input}
\SetKwInOut{Output}{Output}
\Input{$N$, $K$, $T$, $n$, $m$, $s=0$, $u=0$, $v_0 = 0$}
\Output{$d\langle J\rangle/ds \approx (s + u)/N$}

Randomly generate: $x_0$, $Q_0$, $a^{i,j}_0$, $w^i_0$ for all $i,j=1,...,m$\;

\For(\tcp*[h]{main time loop}){$k = 0,...,N-1$}{

\If{$k\geq T$}
{$s = s + DJ_k\cdot v_k$\;
$u = u - J_k\,(u_k + u_{k-1} + ... + u_{k-K+1})$\;}

$S_k = D\varphi_k\,Q_k$\;
QR-factorize $S_k$: $Q_{k+1}\,R_{k+1} = S_k$\;
Find the inverse of $R_{k+1}$\;

\lFor{$i = 1,...,m$, $j= 1,...,i$}{
$\tilde{a}_{k+1}^{i,j} = D^2\varphi_k(q_k^{i}, q_k^{j}) + D\varphi_k\,a_{k}^{i,j}$
}

\lFor{$i = 1,...,m$, $j = 1,...,i$}{
$a_{k+1}^{i,j} = \tilde{a}_{k+1}^{p,q}\,(R^{-1})_{k+1}^{(pi)}\,(R^{-1})_{k+1}^{(qj)}$
}

\For{$i = 1,...,m$}{
\For{$p, q = 1,...,m$}{
$(\partial_{\xi_{k+1}^i}R_{k+1})^{(pq)} = \begin{cases} q_{k+1}^p\cdot a_{k+1}^{p,i}, & \text{if } p = q \\ q_{k+1}^p\cdot a_{k+1}^{q,i} + q_{k+1}^{q}\cdot a_{k+1}^{p,i}, & \text{if } p < q \\
    0, & \text{otherwise}\end{cases}
$\;
}
$g_{k+1}^i = -\text{tr}(\partial_{\xi_{k+1}^i}R_{k+1})$\;
}

$f_k = D\varphi_k\,v_{k} + \chi_{k+1}$\;

\For{$i = 1,...,m$}{
$c_{k+1}^i = q_{k+1}^i\cdot f_k$\;
$\partial_{\xi_k^i}\,f_k = D^2\varphi_k(v_k,q_k^i) + D\varphi_k\,w_k^i + D\partial_s\varphi_k\,q_k^{i}$\;
}

$v_{k+1} = f_k - c_{k+1}^i\,q_{k+1}^i$\;
$\nabla_{\xi_{k+1}}f_k = \nabla_{\xi_{k}}f_k\,R_{k+1}^{-1}$\;

\For{$i,j = 1,...,m$}{
$p_{k+1}^{i,j} = a_{k+1}^{i,j} - q_{k+1}^{l}(\partial_{\xi_{k+1}^j}R_{k+1})^{(li)}$\;
$b_{k+1}^{i,j} = p^{i,j}_{k+1}\cdot f_{k} + q_{k+1}^i\cdot (\nabla_{\xi_{k+1}}f_k)^{:j}$\;
}

\lFor{$i = 1,...,m$}{
$w_{k+1}^{i} = (\nabla_{\xi_{k+1}}f_k)^{:i} - b_{k+1}^{l,i}\,q_{k+1}^l + c_{k+1}^l\,p_{k+1}^{l,i}$
}

Save the scalar: $u_{k+1} = b^{i,i}_{k+1} + c_{k+1}^i\,g_{k+1}^i$\;
Advance the iteration: $x_{k+1} = \varphi(x_k)$\;
Evaluate: $D\varphi_{k+1}$, $D^2\varphi_{k+1}$, $D\partial_s\varphi_{k+1}$, $\chi_{k+2}$ , $J_{k+1}$, $DJ_{k+1}$\;

}
\caption{Space-split sensitivity (S3) algorithm for discrete systems}
\end{algorithm}

Every iteration of the main time loop starts from updating the sums of the stable and unstable integrands, $s$ and $u$, respectively (Lines 3-6). We disregard first $T$ data points to ensure all quantities contributing to the final average are close to their respective true values up to the machine precision. Given all recursions exponentially converge, the value of $T$ is in fact relatively low and $T \ll N$.    

Lines 7-11 are taken from \cite{sliwiak-srb}, as they reflect all the steps necessary to compute the SRB density gradient. Note this code chunk involves advancing $m$ tangent equations (Line 7), QR factorization (Line 8), inverting the $R$ matrix (Line 9), advancing $m^2/2$ second-order tangent equations (Line 10), and double rescaling of $m^2/2$ vectors against the $R^{-1}$ matrix (Line 11). Indeed, the most expensive stage of this chunk is Line 10, which costs $\mathcal{O}(n^3\,m^2)$ due to the presence of the third-order tensor (Hessian of $\varphi$) contracted against two different vectors. This is because for each component of the new $n$-dimensional vector $\tilde{a}$, one must compute and sum up $n^2$ different scalar products. Note also that the rescaling stage (Line 11) involves four nested for-loops, which implies the brute-force vector-by-vector rescaling would require $\mathcal{O}(n\,m^4)$ floating point operations. However, as pointed out in Section \ref{sec:unstable}, this operation can also be completed in a component-by-component fashion. Indeed, the rescaling process involves double contraction against the same matrix. It means that one needs to compute $n$ matrix products $R_{k+1}^T\,\tilde{A}^i\,R_{k+1}$, where $\tilde{A}^i$ denotes an $m\times m$ matrix containing $i$-th components of all $m^2$ vectors $\tilde{a}$. The double matrix-matrix product costs $\mathcal{O}(m^3)$ flops and therefore the total cost of Line 11 is proportional to $\mathcal{O}(n\,m^3)$. The reader is referred to \cite{sliwiak-srb} for a more detailed analysis of the computational complexity of this part of the algorithm. 

Lines 12-17 compute $m$ upper-triangular derivatives of the $R$ matrix. Each component requires evaluating one or two dot products, which implies the cost of executing this chunk is $\mathcal{O}(n\,m^3)$. Here, we automatically obtain the SRB density gradient by evaluating the traces of all $\partial_{\xi}R$. The simplified relation for $g$ (Line 16) is a direct consequence of the measure conservation, which was obtained through parametric differentiation of Eq. \ref{eqn:unstable-measure-convervation} using locally orthogonal coordinates (a complete derivation can be found in \cite{sliwiak-densitygrad}). The leading order of the flop count of the code fragment involving Lines 18-23 is determined by Line 21. This line evaluates $m$ parametric derivatives of $f$ through the Hessian contraction and two other matrix-vector products. Per our discussion above, therefore, Line 21 requires $\mathcal{O}(n^3\,m)$ flops. The same estimate also applies to the algorithm part involving Lines 24-30. Here, the most expensive stage is Line 26, which evaluates $m^2$ vectors $p$. Note for each vector $p$, we compute a matrix-vector product, each requiring $\mathcal{O}(nm)$ algebraic operations. Note also the variable change (Line 24), computation of all scalars $b$ (Line 27), and the update of $w$ (Line 29) cost $\mathcal{O}(n\,m^2)$ each. The final chunk of this algorithm, Lines 31-32, evaluates a collection of nonlinear expressions and thus its complexity depends on the structure of the system.

Since $n\geq m$, the leading term of the total flop count of Algorithm \ref{alg:alg1} (excluding the nonlinear part) is proportional to $\mathcal{O}(n^3\,m^2)$. This estimate reflects the worst-case scenario, for a general chaotic system. Many real-world chaotic systems, however, produced by spatial discretization of partial differential equations (PDEs), have a special structure. Popular discretization schemes, such as the finite element method, generate systems with local dependencies. Each grid point is usually communicated only with neighboring points. This implies both the Jacobian and Hessian of $\varphi$ are sparse arrays and have a banded structure. In such systems, therefore, the total flop count is linear with respect to the dimension of the system $n$. To conclude, the ultimate cost of approximating the sensitivity $d\langle J\rangle/ds$ using Algorithm \ref{alg:alg1} and data from $N$ consecutive states along a typical trajectory is $\mathcal{O}(N\,n^3\,m^2)$. For physical systems, however, this estimate can potentially be reduced to $\mathcal{O}(N\,n\,m^3)$ in the presence of sparsity patterns arising due to local discretization. 

In terms of the storage, the largest arrays are the following: the Hessian $D^2\varphi$, $m^2$ vectors $p$, and $m^2/2$ vectors $a$. They collectively have $n^3$, $m^2\,n$ and $m^2\,n$ components, respectively. In case of PDE-related systems with a sparse structure, the number of components to be stored is linear with respect to $n$. Moreover, in several physical dynamical systems, the dimension of the unstable manifold is significantly smaller than the system's dimension, i.e., $m/n \ll 1$ and $m^2 < n$ \cite{blonigan-phdthesis,ni-jfm, blonigan-ks}. Thus, in case of sparse physical systems, our algorithm requires storing two arrays no larger than $n^2$ entries and a few significantly smaller arrays. Note that in order to compute all required quantities at step $k+1$, we only need information from the previous one, i.e., $k$-th time step. No information from steps $k-2,k-3,...,0$ is required to advance the iteration.

\section{Numerical results}\label{sec:results}


The purpose of this section is to test Algorithm \ref{alg:alg1} using two low-dimensional chaotic maps taken from the literature. In particular, we shall consider the two-dimensional ($n=2$) baker's map $\varphi:[0,2\pi]^2\to[0,2\pi]^2$ \cite{chandramoorthy-s3-new},
\begin{equation}
\label{eqn:results-bakers}
    \begin{split}
        x_{k+1}^{(1)} &= 2x_k^{(1)} + s_{1}\,\sin(x_{k}^{(1)}) +  s_{2}\,\sin(x_{k}^{(1)})\,\sin(2x_{k}^{(2)})/2 \,\text{mod}\,2\pi\\
        x_{k+1}^{(2)} &= x_k^{(2)}/2 + \pi \lfloor x_k/\pi\rfloor  + s_3\,\sin(x_{k}^{(1)})\,\sin(2x_k^{(2)})/2 + s_4\,\sin(2x_k^{(2)})/2 \,\text{mod}\,2\pi
    \end{split},
\end{equation}
and the extended three-dimensional ($n=3$) solenoid map \cite{ni-fast, williams-solenoid},
\begin{equation}
    \label{eqn:results-nimap}
    \begin{split}
        x_{k+1}^{(1)} &= 0.05x_k^{(1)} + 0.1\cos(8x_k^{(2)}) - 0.1\sin(5x_k^{(3)})\\
        x_{k+1}^{(2)} &= 2x_k^{(2)} + s\,(1+x_k^{(1)})\,\sin(8x_k^{(2)})\,\text{mod}\,2\pi \\
        x_{k+1}^{(3)} &= 3x_k^{(3)} + s\,(1+x_k^{(1)})\,\cos(2x_k^{(3)})\,\text{mod}\,2\pi
    \end{split}.
\end{equation}

The baker's map involves a set of four real-valued parameters, $\{s_1,s_2,s_3,s_4\}$. For moderately low parameter values, this map has one positive LE ($m=1$) close to $\log 2$ and one negative LE. Eq. \ref{eqn:results-bakers} is a mathematical representation of the kneading operation, in which a thin dough is stretched by the factor of 2 and then compressed by the same factor. This stretching/compressing process is perturbed in two directions through the sine functions. Baker's maps serve as deterministic models of diffusion processes and are widely used in statistical mechanics \cite{gaspard-bakers}. The second map, in Eq. \ref{eqn:results-nimap}, is parameterized by a single real-valued parameter $s$. It was constructed in \cite{ni-fast} by adding one additional expanding rotation and extra interaction terms between contracting and expanding directions of the Smale-Williams map used in modeling of oscillating circuits \cite{kuznetsov-solenoid}. If $s$ is moderately low, this map has two positive LEs ($m=2$), with values close to $\log 2$ and $\log 3$, and a negative one. Therefore, unstable manifolds are geometrically represented by smooth curves immersed in $\mathbb{R}^2$ (baker's map) and surfaces immersed in $\mathbb{R}^3$ (solenoid map).  

\begin{figure}
    \centering
    \begin{minipage}{0.49\textwidth}
        \includegraphics[width = \textwidth]{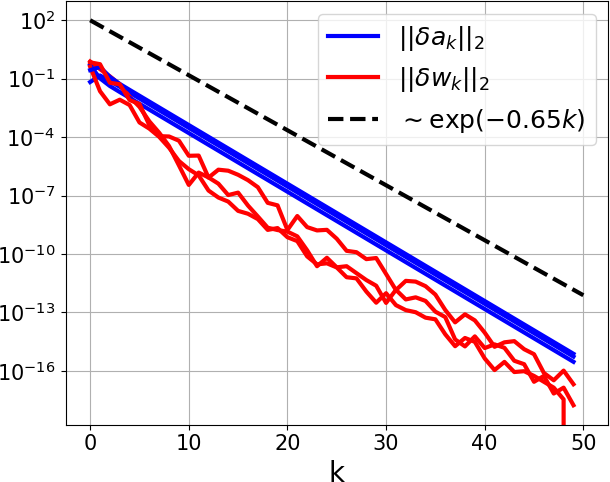}
    \end{minipage}
    \begin{minipage}{0.49\textwidth}
        \includegraphics[width = \textwidth]{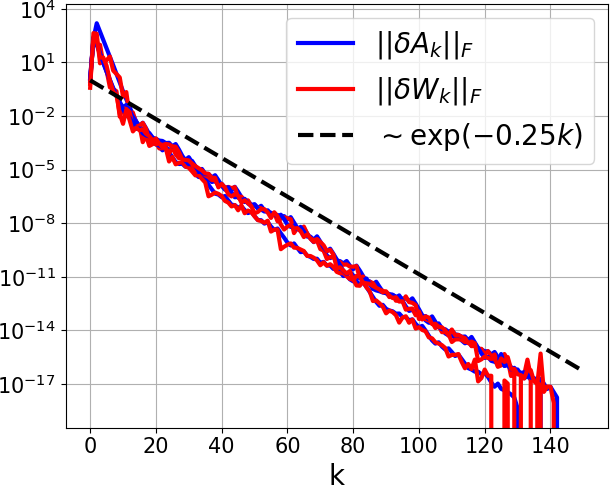}
    \end{minipage}
    \caption{Convergence of the iterative formulas for $a$ and $w$ on the semi-logarithmic scale. We compute vector 2-norms and matrix Frobenius norms, respectively for the baker's map (left) and solenoid map (right), of the difference arrays as a function of the time step $k$. In both cases, the parameter values are randomly chosen from the uniform distribution $[0,0.2]$. The black dashed lines are reference lines proportional to an exponential function. Each line of the same color corresponds to a different trajectory.}
    \label{fig:results-conv-aw}
\end{figure}

Figure \ref{fig:results-conv-aw} presents convergence plots of the recursive formulas for $a$ and $w$, which are key ingredients of Algorithm \ref{alg:alg1}. We perform this test by randomly choosing two different initial conditions $w_{0}^i$ and $a_{0}^{i,j}$, $i,j = 1,...,m$ and tracing their respective difference vectors (as defined in Eq. \ref{eqn:unstable-a-diff} and Eq. \ref{eqn:unstable-w-conv}) as a function of time $k$. We also randomly choose parameter values and repeat this test three times (i.e., by following three different trajectories). 

For both the baker's map and solenoid map, the quantities obtained through the recursions derived in Section \ref{sec:unstable} exponentially converge which confirms our analytical predictions. We observe that the rate of exponential convergence may vary from system to system. In case of the solenoid map, we notice a significant peak right after the beginning of the recursion. This is a consequence of the randomly chosen initial condition $x_0$ that is likely to be located beyond the attractor, given its complex geometry \cite{ni-fast}.

\begin{figure}
    \centering
    \begin{minipage}{0.49\textwidth}
        \includegraphics[width = \textwidth]{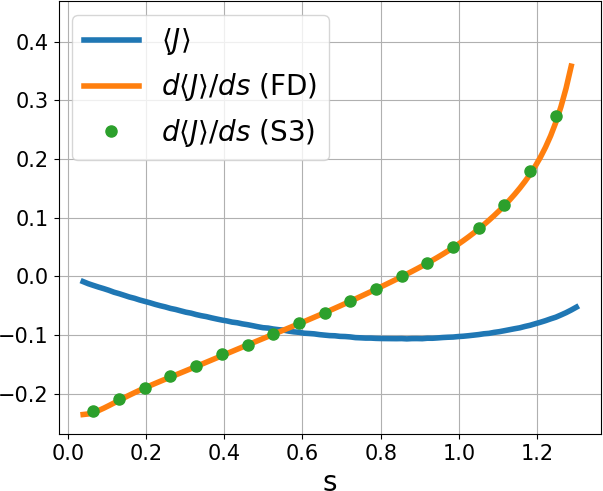}
    \end{minipage}
    \begin{minipage}{0.49\textwidth}
        \includegraphics[width = \textwidth]{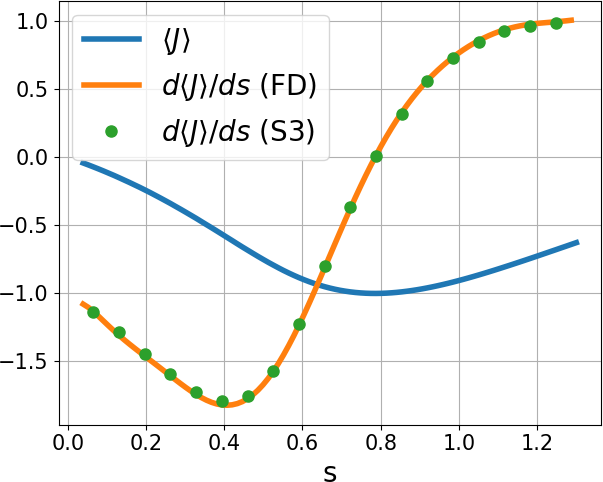}
    \end{minipage}
    \caption{Sensitivity approximation of the baker's map (Eq. \ref{eqn:results-bakers}) using Algorithm \ref{alg:alg1}. In this particular simulation, we set $s=s_1=s_2$, $s_3=s_4=0$, $N=10^6$, $K=11$, $T=100$, and choose $J = \cos(4x^{(2)})$. The long-time averages $\langle J\rangle$ were generated through the 11-degree polynomial fit of an evenly spaced data set, where each data point was computed with $N=10^8$ samples. The solution reference line (FD) approximates the sensitivity and was obtained by the central finite difference applied to the curve fit.}
    \label{fig:results-2D-sensitivity}
\end{figure}
\begin{figure}
    \centering
    \begin{minipage}{0.49\textwidth}
        \includegraphics[width = \textwidth]{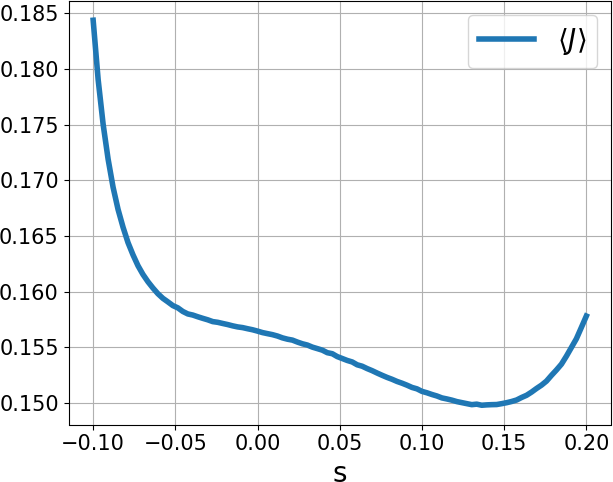}
    \end{minipage}
    \begin{minipage}{0.49\textwidth}
        \includegraphics[width = \textwidth]{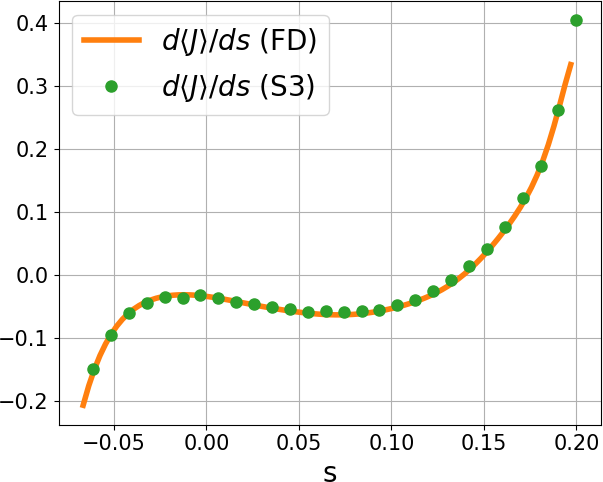}
    \end{minipage}
    \caption{Sensitivity approximation of the solenoid map (Eq. \ref{eqn:results-nimap}) using Algorithm \ref{alg:alg1}. In this particular simulation, we set $N=10^7$, $K=11$, $T=100$, and choose $J = \sin(x^{(2)})\,\cos(4x^{(2)})\,x^{(3)}$. The long-time average $\langle J\rangle$ and its finite difference approximation of the parametric derivative (FD) were generated in the same fashion as their counterparts in Figure \ref{fig:results-2D-sensitivity}.}
    \label{fig:results-3D-sensitivity}
\end{figure}

Given the convergence test results, we set $T=100$ (cut-off threshold) in the space-split algorithm to approximate sensitivities of both the maps. Figures \ref{fig:results-2D-sensitivity}-\ref{fig:results-3D-sensitivity} show the computed approximations generated using Algorithm \ref{alg:alg1} (S3). We validate all S3 outputs by comparing them against their respective finite difference sensitivity approximations (FD). The latter are obtained by computing long-time averages of a chosen objective function, $\langle J\rangle$, at $s\pm \delta s$ and applying the central finite difference scheme with $\delta s = 0.01$. The solid blue and orange lines in Figures \ref{fig:results-2D-sensitivity}-\ref{fig:results-3D-sensitivity} are in fact polynomial curve fits, which were computed only for demonstration purposes. More technical details are included in the captions of these two figures. The key message of these simulations is that Algorithm \ref{alg:alg1} generates accurate sensitivity estimates in discrete systems with an arbitrary number of Lyapunov exponents $m$. 

\begin{figure}
    \centering
    \begin{minipage}{0.325\textwidth}
    \includegraphics[width = 0.94\textwidth]{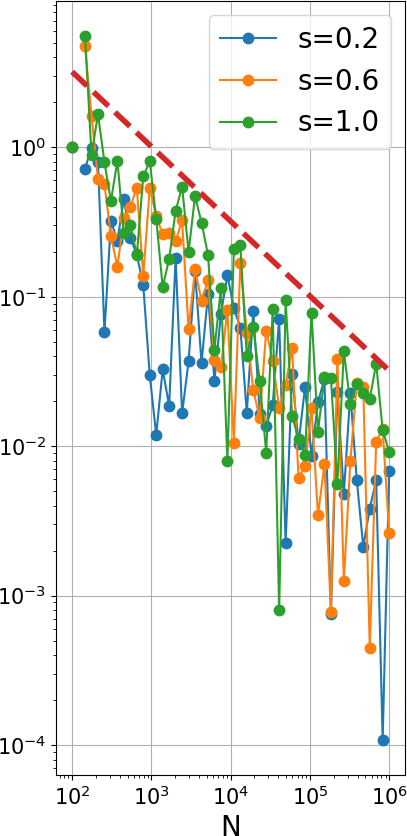}
    \end{minipage}
    \begin{minipage}{0.325\textwidth}
    \includegraphics[width = 1\textwidth]{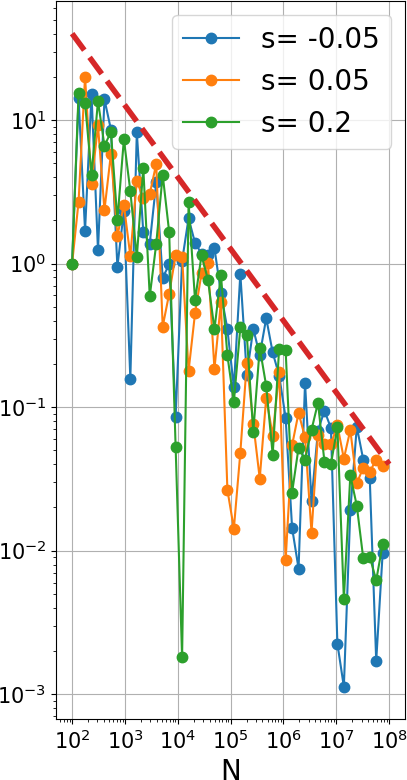}
    \end{minipage}
    \begin{minipage}{0.325\textwidth}
    \includegraphics[width = 0.94\textwidth]{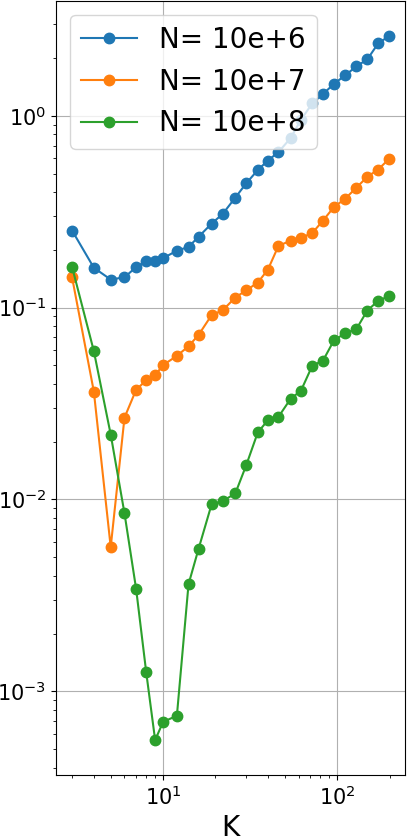}
    \end{minipage}

    \caption{Relation between the relative error of the S3 approximation and $N$ or $K$. The left and middle plots correspond to the baker's (with $s=s_1=s_2$, $s_3=s_4=0$) and solenoid map, respectively, while the right one represents the solenoid map at $s=-0.05$. All relative errors were computed with respect to the finite difference sensitivity approximations generated using $N=10^{10}$ data points. The dashed red lines represent reference lines with the slope $-1/2$ in the loglog scaling.}
    \label{fig:results-conv-NK}
\end{figure}

While we already know the iterative formulas of the space-split algorithm converge exponentially, per our analysis and numerical evidence shown in Figure \ref{fig:results-conv-aw}, the overall accuracy of our method depends on the amount of data used in the ergodic-averaging of the computed time series. In particular, we have freedom to tune the accuracy by modifying the values of $N$ (trajectory length) and $K$ (infinite series truncation number). Figure \ref{fig:results-conv-NK} shows the dependence of the relative error of the S3 approximation on both the parameters.

To produce the left and central plot of Figure \ref{fig:results-conv-NK}, we fixed $K=20$ and computed the relative errors for both the baker's and solenoid map, respectively, at different parameter values. These two plots clearly confirm our algorithm behaves as a typical Monte Carlo scheme, i.e., the error scales as $\mathcal{O}(1/\sqrt{N})$. Notice also we respectively need $N=\mathcal{O}(10^6)$ and $N=\mathcal{O}(10^8)$ data points to secure the relative error $\mathcal{O}(10^{-2})$, which means the constant $C_1$ from Ineq. \ref{eqn:unstable-approximation-error} may significantly change from system to system. The right-hand side plot of Figure \ref{fig:results-conv-NK} indicates the relation of the error and $K$ (for a fixed $N$) is non-monotonic, which is consistent with the rigorous estimate of Ineq. \ref{eqn:unstable-approximation-error}. Indeed, we observe the interaction of the linear and exponential functions of that inequality. If $N$ is sufficiently large but finite, we always observe sudden (exponential) decay of the error for all $0\leq K\leq K^*$ and then, for all $K\geq K^*$, proportional increase of the error. Certainly, $K^*$ depends on the model itself, but also on the value of $N$, as shown in the plot. In practice, one can store several approximations of the unstable contribution (each corresponding to a different $K$; ideally, the chosen values of $K$ are logarithmically separated) and choose the one that significantly breaks the monotonicity. The cost of computing multiple approximations of the unstable contribution is negligible compared to the total cost of Algorithm \ref{alg:alg1}.      

\section{Conclusions}\label{sec:conclusions}

It is generally difficult to accurately estimate sensitivities of chaotic dynamical systems. Due to the butterfly effect, the direct simulation of solution perturbations is impractical. Several numerical methods have been proposed to compute the sensitivity of chaos, but most of them suffer from at least one of the following common problems: exploding tangent solutions, unphysicality of shadowing trajectories, huge computational cost and storage requirements, complicated generalization. 

Our new method for sensitivity analysis derives from Ruelle's rigorous linear response theory, which is regularized based the concept of perturbation space-splitting and partial integration along unstable manifolds \cite{chandramoorthy-s3-new}. Through the intuitive measure-based parameterization of the unstable subspace \cite{sliwiak-densitygrad,sliwiak-srb} and chain rule on smooth manifolds, we systematically derive a set of iterative (trajectory-following) formulas for different quantities arising in the regularization of Ruelle's formula, and show their exponential convergence. Similarly to the majority of methods that stem from the linear response theory, our method is formulated as a typical Monte Carlo scheme, which rigorously converges to the true solution as $\mathcal{O}(1/\sqrt{N})$, where $N$ is the trajectory length. The following list summarizes the main advantages of the space-split approach:
\begin{itemize}
\item Immunity to the ergodicity-breaking/unphysicality errors (common in some approximative methods) and the omnipresent butterfly effect,
\item Generalizability to $n$-dimensional systems, $n\in\mathbb{Z}^+$; the algorithm we propose is ready-to-use for discrete systems with an arbitrary number of positive Lyapuonv exponents $m$, 
\item Provable convergence for uniformly hyperbolic systems,
\item Translatable to memory-efficient as-we-go Monte Carlo algorithms.
\end{itemize}
The major consequence of partial Lebesgue integration is the computation of directional derivatives of an ergodic measure, describing the system's statistical behavior. Recursive computation of this quantity, known as the SRB density gradient, requires solving a collection first- and second-order tangent equations. Indeed, this is the actual price of the regularization of Ruelle's formula. From the algorithmic perspective, therefore, we must perform a series of algebraic operations involving third-order tensors. We estimate the total flop count is $\mathcal{O}(n^3\,m^2)$, which can be reduced to $\mathcal{O}(n\,m^3)$ for PDE-related systems with local dependencies and sparse structures.

While this paper solely focuses on discrete systems, the proposed algorithm can naturally be extended to continuous-in-time (ODE) systems. The perturbation vector splitting would need to incorporate the one-dimensional neutral subspace that is aligned with the flow. This requires the derivation of extra recursive formulas for certain new quantities and addition of the neutral contribution to the sensitivity approximation. Given the one-dimensionality of the extra subspace, the leading term of the algorithm's total flop count is expected to remain the same.

\section*{Supplementary Material}

To facilitate the reproduction of all reported results, the authors attach the Python code used to generate Figures \ref{fig:results-conv-aw}-\ref{fig:results-conv-NK}. Inside the main folder, the reader will find a ``README" file containing a description of all Python scripts. 

\section*{Acknowledgments}
This work was funded by Air Force Office of Scientific Research Grant No. FA8650-19-C-2207 and U.S. Department of Energy Grant No. DE-FOA-0002068-0018. The authors also acknowledge the MIT SuperCloud and Lincoln Laboratory Supercomputing Center for providing HPC resources that have contributed to the research results reported within this paper.

\section*{Conflict of interest}
The authors declare that they have no conflict of interest.

\bibliographystyle{elsarticle-num-names}
\bibliography{references.bib}

\end{document}